\documentclass[12pt,oneside,leqno]{amsproc}
\usepackage[utf8]{inputenc}
\usepackage[russian]{babel}
\usepackage{amssymb,amsmath,amsthm}
\renewcommand{\subsection}{\refstepcounter{subsection}%
\par\bigskip\noindent\textbf{\upshape\thesubsection. }}
\renewcommand{\subsubsection}{\refstepcounter{subsubsection}%
\par\medskip\noindent\textbf{\upshape\thesubsubsection.  }}
\renewcommand{\paragraph}{\refstepcounter{paragraph}%
\par\smallskip\noindent\textbf{\upshape\theparagraph. }}

\numberwithin{equation}{subsection}

\renewcommand{\thesubsection}{\arabic{subsection}}
\renewcommand{\thesubsubsection}{\arabic{subsection}.\arabic{subsubsection}}

\makeatletter
\renewcommand{\@makefnmark}{\hbox{\small\mathsurround=0cm%
${}\hspace{0.04cm}{}^{\@thefnmark)}$}}
\renewcommand{\@makefntext}[1]{\parindent=1em\noindent\hbox to 1.8em{%
\hss${}^{\@thefnmark)}$}\,#1}
\makeatother

\title{Осцилляционный метод в задаче о спектре дифференциального оператора
четвёртого порядка с самоподобным весом}
\author{А.~А.~Владимиров\footnote{Работа поддержана РФФИ, грант \No~10-01-00423.}}
\begin{document}
\renewcommand{\proofname}{{\upshape Д\,о\,к\,а\,з\,а\,т\,е\,л\,ь\,с\,т\,в\,о.}}
\begin{abstract}
Рассматриваются самосопряжённые граничные задачи для дифференциального выражения
\[
	y^{(4)}-\lambda\rho y=0,
\]
где вес \(\rho\in W_2^{-1}[0,1]\) представляет собой обобщённую производную
самоподобной функции канторовского типа. На основе изучения осцилляционных
свойств собственных функций уточняются характеристики известных спектральных
асимптотик таких задач.
\end{abstract}
\begin{flushleft}
\normalsize УДК~517.984
\end{flushleft}
\maketitle
\markboth{}{}

\section{Введение}\label{par:1}
\subsection\label{pt:1.1}
Целью настоящей статьи является применение разработанного в \cite{VSh:2011}
осцилляционного метода исследования спектральных асимптотик задач с самоподобными
весами к случаю самосопряжённой граничной задачи
\begin{gather}\label{ekvac}
	y^{(4)}-\lambda\rho y=0,\\ \label{grusl} (U-1)y^{\vee}+i(U+1)y^{\wedge}=0,
\end{gather}
где \(\rho\in W_2^{-1}[0,1]\) "--- неотрицательная обобщённая весовая функция,
\(U\in\mathbb C^{4\times 4}\) "--- унитарная матрица граничных условий,
а \(y^{\wedge}\) и \(y^{\vee}\) "--- числовые векторы
\[
	y^{\wedge}\rightleftharpoons\begin{pmatrix}y^{[0]}(0)&y^{[1]}(0)&
		y^{[0]}(1)&y^{[1]}(1)\end{pmatrix}^T,\qquad
	y^{\vee}\rightleftharpoons\begin{pmatrix}y^{[3]}(0)&y^{[2]}(0)&
		-y^{[3]}(1)&-y^{[2]}(1)\end{pmatrix}^T.
\]
Через \(y^{[k]}\), где \(k\in\{0,\ldots,3\}\), здесь обозначены стандартные
\cite[\S~15]{Na:1969}, \cite[(7.46)]{RH:2001} квазипроизводные
\(y^{[0]}\rightleftharpoons y\), \(y^{[1]}\rightleftharpoons y'\),
\(y^{[2]}\rightleftharpoons y''\) и \(y^{[3]}\rightleftharpoons -y'''\).
Содержание работы \cite{VSh:2011} будет далее предполагаться известным.

\subsection\label{pt:grusl}
Граничные задачи \ref{pt:1.1}\,\eqref{ekvac}, \ref{pt:1.1}\,\eqref{grusl} будут
далее рассматриваться не в максимальной общности. А именно, соотношения
\ref{pt:1.1}\,\eqref{grusl} мы намерены предполагать допускающими запись в виде
\begin{multline}\label{grusl1}
	y^{[2]}(0)+\alpha y^{[0]}(0)-\beta y^{[1]}(0)=\beta y^{[3]}(0)+
		\alpha y^{[2]}(0)=\\ =y^{[2]}(1)+\alpha y^{[0]}(1)+
		\beta y^{[1]}(1)=\beta y^{[3]}(1)-\alpha y^{[2]}(1)=0,
\end{multline}
где \(\alpha\geqslant 0\), \(\beta>0\). Кроме того, функция \(\rho\) будет
обычно предполагаться обобщённой производной неубывающей функции \(P\in C[0,1]\)
\emph{канторовского типа самоподобия}. Это означает \cite[\S~2]{VSh:2011}
выполнение равенств \(P(0)=0\) и \(P(1)=1\), а также существование натурального
числа \(\varkappa>1\) и пары вещественных чисел \(a\in (0,1/\varkappa)\),
\(b\rightleftharpoons (1-\varkappa a)/(\varkappa-1)\) со следующими свойствами:
\begin{enumerate}
\item Независимо от выбора индекса \(k\in\{0,\ldots,\varkappa-1\}\) функция
\(P_k\in C[0,1]\) вида
\[
	P_k(x)\rightleftharpoons\varkappa P(k[a+b]+ax)
\]
совпадает с функцией \(P\) с точностью до аддитивной постоянной.
\item Независимо от выбора индекса \(k\in\{1,\ldots,\varkappa-1\}\) функция \(P\)
постоянна на интервале \((k[a+b]-b,k[a+b])\).
\end{enumerate}
Некоторые факты о распределении спектра задач рассматриваемого типа могут быть
найдены в работе \cite[\S~3]{Naz:2004}.

\subsection\label{pt:1.2}
Формальной задаче \ref{pt:1.1}\,\eqref{ekvac}, \ref{pt:grusl}\,\eqref{grusl1}
обычным образом \cite{V:2004} сопоставляется линейный пучок \(T:\mathbb C\to
\mathcal B(W_2^2[0,1],W_2^{-2}[0,1])\) операторов вида
\begin{equation}\label{kvf}
	\langle T(\lambda)y,y\rangle\equiv\int\limits_0^1 |y''|^2\,dx+
		\dfrac{|\alpha y(0)-\beta y'(0)|^2+|\alpha y(1)+\beta y'(1)|^2}{%
		\beta}-\lambda\langle\rho,|y|^2\rangle.
\end{equation}
Интегрированием по частям \cite[Лемма~2]{V:2004} легко устанавливается, что пара
\(\{\lambda,y\}\) из числа \(\lambda\in\mathbb C\) и нетривиальной функции
\(y\in W_2^2[0,1]\) является собственной парой пучка \(T\) в том и только том
случае, когда функции \(y''\) и \(y'''-\lambda Py\) непрерывно дифференцируемы
и удовлетворяют уравнению
\begin{equation}\label{ekvac1}
	[y'''-\lambda Py]'+\lambda Py'=0
\end{equation}
совместно с понимаемыми в обычном смысле граничными условиями \ref{pt:grusl}\,%
\eqref{grusl1}. Это наблюдение постоянно будет использоваться нами в дальнейшем.

\subsection
Статья имеет следующую структуру. В \ref{par:2} излагаются сведения об осцилляции
собственных функций задач рассматриваемого типа. Они являются достаточно
стандартными \cite{GK:1950}, \cite{LN:1958} и не претендуют в полной мере
на научную новизну. В \ref{par:3} рассматривается явление спектральной
периодичности и вытекающие из него свойства спектральных асимптотик, а также
приводятся иллюстрирующие полученные теоретические результаты данные численных
экспериментов.


\section{Осцилляция собственных функций}\label{par:2}
\subsection
Имеют место следующие два факта:

\subsubsection\label{prop:1.1}
{\itshape Пусть \(\lambda>0\), а \(y\in C^3[0,1]\) есть нетривиальное решение
уравнения \ref{par:1}.\ref{pt:1.2}\,\eqref{ekvac1}, удовлетворяющее при некотором
\(a\in [0,1)\) неравенствам
\[
	y(a)\geqslant 0,\qquad y'(a)\geqslant 0,\qquad y''(a)\geqslant 0,\qquad
	y'''(a)\geqslant 0.
\]
Тогда выполняются также неравенства
\[
	y(1)>0,\qquad y'(1)>0,\qquad y''(1)>0,\qquad y'''(1)>0.
\]
}

\begin{proof}
Зафиксируем последовательность \(\{P_n\}_{n=0}^{\infty}\) равномерно стремящихся
к функции \(P\) функций класса \(C^1[0,1]\), имеющих равномерно положительные
производные и удовлетворяющих равенствам \(P_n(a)=P(a)\). Зафиксируем также
последовательность \(\{y_n\}_{n=0}^{\infty}\) решений начальных задач
\begin{gather}\label{eq:ekv}
	[y'''_n-\lambda P_ny_n]'+\lambda P_ny'_n=0,\\ \notag
	y_n^{(k)}(a)=y^{(k)}(a),\qquad k\in\{0,\ldots,3\}.
\end{gather}
Стандартными методами теории линейных дифференциальных уравнений для вектор-функций
\cite[\S~16]{Na:1969} легко устанавливается факт равномерной на отрезке \([0,1]\)
сходимости последовательностей \(\{y_n\}_{n=0}^{\infty}\), \(\{y'_n\}_{n=0}^{%
\infty}\), \(\{y''_n\}_{n=0}^{\infty}\) и \(\{y'''_n-\lambda P_ny_n\}_{n=0}^{%
\infty}\) к функциям \(y\), \(y'\), \(y''\) и \(y'''-\lambda Py\), соответственно.

Согласно \cite[Lemma~2.1]{LN:1958}, каждая из функций \(y_n\), \(y'_n\), \(y''_n\)
и \(y'''_n\) строго положительна на полуинтервале \((a,1]\). Объединяя этот факт
с уравнениями \eqref{eq:ekv}, устанавливаем справедливость оценок
\[
	(\forall x\in (a,1))\qquad y'''_n(1)\geqslant\lambda[P_n(1)-P_n(x)]y_n(x).
\]
Посредством предельного перехода теперь немедленно устанавливается факт
неотрицательности на полуинтервале \((a,1]\) каждой из функций \(y\), \(y'\),
\(y'\) и \(y'''\), а также справедливость оценок
\begin{equation}\label{eq:taks}
	(\forall x\in (a,1))\qquad y'''(1)\geqslant\lambda[P(1)-P(x)]y(x).
\end{equation}
При этом, ввиду нетривиальности функции \(y\), заведомо найдётся величина
\(\gamma>0\) со свойством
\begin{equation}\label{eq:taks1}
	(\forall x\in (a,1])\qquad y(x)\geqslant\gamma\cdot(x-a)^3.
\end{equation}
Объединяя оценки \eqref{eq:taks} и \eqref{eq:taks1} с фактом непостоянности
функции \(P\) в любой левой окрестности точки \(1\), убеждаемся в выполнении
неравенства \(y'''(1)>0\), а тогда и прочих требуемых неравенств.
\end{proof}

\subsubsection\label{prop:1.2}
{\itshape Пусть \(\lambda>0\), а \(y\in C^3[0,1]\) есть нетривиальное решение
уравнения \ref{par:1}.\ref{pt:1.2}\,\eqref{ekvac1}, удовлетворяющее при некотором
\(a\in (0,1]\) неравенствам
\[
	y(a)\geqslant 0,\qquad y'(a)\leqslant 0,\qquad y''(a)\geqslant 0,\qquad
	y'''(a)\leqslant 0.
\]
Тогда выполняются также неравенства
\[
	y(0)>0,\qquad y'(0)<0,\qquad y''(0)>0,\qquad y'''(0)<0.
\]
}

Утверждение \ref{prop:1.2} доказывается полностью аналогично утверждению
\ref{prop:1.1}.

\subsection
Имеют место следующие два факта:

\subsubsection\label{prop:2.4}
{\itshape Спектр пучка \(T\) составлен последовательностью \(\{\lambda_n
\}_{n=0}^{\infty}\) неотрицательных "--- а в случае \(\alpha>0\) даже строго
положительных "--- простых собственных значений. Независимо от выбора индекса
\(n\in\mathbb N\) отвечающая собственному значению \(\lambda_n\) собственная
функция \(y_n\) имеет только простые нули и удовлетворяет условиям \(y_n(0)\neq 0\)
и \(y_n(1)\neq 0\).
}

\begin{proof}
Из выражения \ref{par:1}.\ref{pt:1.2}\,\eqref{kvf} квадратичной формы оператора
\(T(0)\) немедленно вытекает, что ядро этого оператора образовано линейными
функциями, удовлетворяющими равенствам
\[
	\alpha y(0)-\beta y'(0)=\alpha y(1)+\beta y'(1)=0.
\]
В случае \(\alpha>0\) единственной такой функцией является тождественно нулевая.
Соответственно, в этом случае все собственные значения пучка \(T\) строго
положительны. В случае \(\alpha=0\) такие функции образуют одномерное
подпространство постоянных функций.

Пусть некоторое собственное значение \(\lambda>0\) пучка \(T\) обладает
собственной функцией \(y\) со свойством \(y(0)=0\). Тогда без ограничения общности
рассмотрения можно считать, что функция \(y\) вещественнозначна, а знаки величин
\(y'(0)\), \(y''(0)\) и \(y'''(0)\) совпадают [\ref{par:1}.\ref{pt:grusl}\,%
\eqref{grusl1}]. Однако это влечёт противоречащее граничным условиям
\ref{par:1}.\ref{pt:grusl}\,\eqref{grusl1} совпадение знаков величин \(y''(1)\neq
0\) и \(y'''(1)\neq 0\) [\ref{prop:1.1}].

Пусть некоторое собственное значение \(\lambda>0\) пучка \(T\) обладает
собственной функцией \(y\) со свойством \(y(1)=0\). Тогда без ограничения общности
рассмотрения можно считать, что функция \(y\) вещественнозначна, а знаки величин
\(-y'(1)\), \(y''(1)\) и \(-y'''(1)\) совпадают [\ref{par:1}.\ref{pt:grusl}\,%
\eqref{grusl1}]. Однако это влечёт противоречащее граничным условиям
\ref{par:1}.\ref{pt:grusl}\,\eqref{grusl1} совпадение знаков величин \(y''(0)\neq
0\) и \(-y'''(0)\neq 0\) [\ref{prop:1.2}].

Пусть некоторое собственное значение \(\lambda>0\) пучка \(T\) является кратным.
Тогда для него найдётся собственная функция \(y\) со свойством \(y(0)=0\),
что противоречит сказанному ранее.

Наконец, пусть для некоторого собственного значения \(\lambda>0\) пучка \(T\)
существует кратный нуль \(a\in (0,1)\) соответствующей собственной функции \(y\).
Тогда знаки величин \(y''(a)\) и \(y'''(a)\) являются либо совпадающими, либо
различными. Первый случай означает противоречащее граничным условиям
\ref{par:1}.\ref{pt:grusl}\,\eqref{grusl1} совпадение знаков величин
\(y''(1)\neq 0\) и \(y'''(1)\neq 0\) [\ref{prop:1.1}]. Второй случай означает
противоречащее граничным условиям \ref{par:1}.\ref{pt:grusl}\,\eqref{grusl1}
совпадение знаков величин \(y''(0)\neq 0\) и \(-y'''(0)\neq 0\) [\ref{prop:1.2}].
\end{proof}

\subsubsection\label{prop:2.5}
{\itshape В случае \(\alpha>0\) оператор \([T(0)]^{-1}T':W_2^2[0,1]\to W_2^2[0,1]\)
не увеличивает числа перемен знака никакой вещественнозначной функции.
}

\begin{proof}
Ввиду непрерывной в смысле равномерной операторной топологии зависимости оператора
\([T(0)]^{-1}T'\) от выбора весовой функции \(\rho\in W_2^{-1}[0,1]\), достаточно
рассмотреть случай, когда функция \(\rho\) непрерывна и равномерно положительна.
Иначе говоря, достаточно установить, что независимо от выбора натуральных чисел
\(n>0\) и \(m\) наличие у удовлетворяющей граничным условиям \ref{par:1}.%
\ref{pt:grusl}\,\eqref{grusl1} вещественнозначной функции \(y\in C^4[0,1]\)
не менее \(n+m\) перемен знака влечёт наличие не менее \(n\) перемен знака
у функции \(y^{(4)}\). В случае \(m\geqslant 4\) этот факт немедленно вытекает
из теоремы Лагранжа о среднем значении. Общий случай будет рассмотрен нами
на основе метода арифметической индукции.

Итак, пусть известно, что наличие у произвольной удовлетворяющей граничным условиям
\ref{par:1}.\ref{pt:grusl}\,\eqref{grusl1} вещественнозначной функции
\(y\in C^4[0,1]\) не менее \(n+m+1\) перемен знака заведомо влечёт наличие не менее
\(n\) перемен знака у функции \(y^{(4)}\). Пусть также некоторая вещественнозначная
функция \(y\in C^4[0,1]\) удовлетворяет граничным условиям \ref{par:1}.%
\ref{pt:grusl}\,\eqref{grusl1}, и пусть найдутся \(n+m+1\) упорядоченных
по возрастанию точек
\[
	0<\xi_{0,1}<\ldots<\xi_{0,n+m+1}<1
\]
со свойствами \(y(\xi_{0,k})\cdot y(\xi_{0,k+1})<0\), где \(k\in\{1,\ldots,n+m\}\).

Согласно теореме Лагранжа, найдутся \(n+m\) точек \(\xi_{1,k}\in (\xi_{0,k},
\xi_{0,k+1})\) со свойствами \(y'(\xi_{1,k})\cdot y(\xi_{0,k})<0\). При этом либо
найдётся точка \(\xi\in (0,\xi_{0,1})\) со свойством \(y(\xi)\cdot
y(\xi_{0,1})<0\), либо \(y'(0)\cdot y'(\xi_{1,1})<|y'(\xi_{1,1})|^2\), либо
\(y''(0)\cdot y'(\xi_{1,1})>0\). Обоснование указанной альтернативы использует
фигурирующее среди граничных условий \ref{par:1}.\ref{pt:grusl}\,\eqref{grusl1}
выражение величины \(y''(0)\) через \(y(0)\) и \(y'(0)\). В первом случае функция
\(y\) имеет не менее \(n+m+1\) перемен знака, что, по индуктивному предположению,
означает наличие не менее \(n\) перемен знака у функции \(y^{(4)}\). Во втором
и третьем случаях найдётся точка \(\xi_{2,1}\in (0,\xi_{1,1})\) со свойством
\(y''(\xi_{2,1})\cdot y'(\xi_{1,1})>0\). Аналогичным образом, либо найдётся точка
\(\xi\in (\xi_{0,n+m+1},1)\) со свойством \(y(\xi)\cdot y(\xi_{0,n+m+1})<0\),
либо \(y'(1)\cdot y'(\xi_{1,n+m})<|y'(\xi_{1,n+m})|^2\), либо \(y''(1)\cdot
y'(\xi_{1,n+m})<0\). В первом случае функция \(y\) имеет не менее \(n+m+1\) перемен
знака. Во втором и третьем случаях найдётся точка \(\xi_{2,n+m+1}\in
(\xi_{1,n+m},1)\) со свойством \(y''(\xi_{2,n+m+1})\cdot y'(\xi_{1,n+m})<0\).

Объединяя сказанное, получаем, что либо функция \(y^{(4)}\) имеет не менее \(n\)
перемен знака, либо найдутся \(n+m+1\) упорядоченных по возрастанию точек
\[
	0<\xi_{2,1}<\ldots<\xi_{2,n+m+1}<1
\]
со свойствами \(y''(\xi_{2,k})\cdot y''(\xi_{2,k+1})<0\).

Далее, согласно граничным условиям \ref{par:1}.\ref{pt:grusl}\,\eqref{grusl1},
выполняется либо неравенство \(y''(0)\cdot y''(\xi_{2,1})<|y''(\xi_{2,1})|^2\),
либо неравенство \(y'''(0)\cdot y''(\xi_{2,1})>0\). В обоих случаях найдётся точка
\(\xi_{3,0}\in (0,\xi_{2,1})\) со свойством \(y'''(\xi_{3,0})\cdot y''(\xi_{2,1})
>0\). Аналогичным образом, выполняется либо неравенство \(y''(1)\cdot
y''(\xi_{2,n+m+1})<|y''(\xi_{2,n+m+1})|^2\), либо неравенство \(y'''(1)\cdot
y''(\xi_{2,n+m+1})<0\). В обоих случаях найдётся точка \(\xi_{3,n+m+1}\in
(\xi_{2,n+m+1},1)\) со свойством \(y'''(\xi_{3,n+m+1})\cdot y''(\xi_{2,n+m+1})<0\).
Тем самым, функция \(y'''\) имеет не менее \(n+m+1\) перемен знака, что, согласно
теореме Лагранжа, означает наличие не менее \(n+m\geqslant n\) перемен знака
у функции \(y^{(4)}\).
\end{proof}

\subsection
Имеют место следующие два факта:

\subsubsection\label{prop:2.10}
{\itshape Пусть вещественнозначная функция \(f\in W_2^2[0,1]\) удовлетворяет
неравенствам \(f(0)\neq 0\), \(f(1)\neq 0\) и имеет на интервале \((0,1)\) ровно
\(n\), причём простых, нулей. Тогда существует величина \(\varepsilon>0\),
для которой любая вещественнозначная функция \(y\in W_2^2[0,1]\) со свойством
\(\|y-f\|_{W_2^2[0,1]}<\varepsilon\) также имеет на интервале \((0,1)\) ровно
\(n\) простых нулей.
}

\bigskip
Это утверждение тривиальным образом вытекает из факта непрерывности естественного
вложения \(W_2^2[0,1]\hookrightarrow C^1[0,1]\).

\subsubsection\label{prop:2.12}
{\itshape Пусть \(\{\lambda_n\}_{n=0}^{\infty}\) "--- последовательность
занумерованных в порядке возрастания собственных значений пучка \(T\). Тогда
независимо от выбора индекса \(n\in\mathbb N\) отвечающая собственному значению
\(\lambda_n\) собственная функция \(y_n\) имеет в точности \(n\) нулей
на интервале \((0,1)\).
}

\enlargethispage{2\baselineskip}
\begin{proof}
Рассмотрим сначала случай \(\alpha>0\). Заметим, что с каждой вещественнозначной
функцией вида
\begin{equation}\label{eq:100}
	f=\sum\limits_{k=0}^{n} c_ky_k
\end{equation}
можно связать функциональную последовательность \(\{f_m\}_{m=0}^{\infty}\) вида
\[
	f_m\rightleftharpoons\sum\limits_{k=0}^{n} c_k\lambda_k^{m}
		\lambda_n^{-m} y_k.
\]
Пределом этой последовательности в пространстве \(W_2^2[0,1]\) является функция
\(c_ny_n\). Соответственно [\ref{prop:2.5}, \ref{prop:2.4}, \ref{prop:2.10}],
при \(c_n\neq 0\) число знакоперемен функции \(f\) минорирует число нулей функции
\(y_n\). Однако, ввиду линейной независимости семейства собственных функций пучка
\(T\), заведомо найдётся функция вида \eqref{eq:100}, удовлетворяющая условию
\(c_n\neq 0\) и имеющая не менее \(n\) перемен знака на интервале \((0,1)\).
Тем самым, функция \(y_n\) имеет не менее \(n\) нулей.

Далее, зафиксируем нетривиальный вещественнозначный многочлен \(Q\) не превышающей
\(n\) степени, принадлежащий инвариантному подпространству оператора \([T(0)]^{-1}
T'\), которое отвечает дополнительной к набору \(\{\lambda_k\}_{k=0}^{n-1}\) части
спектра. Ввиду бесконечности носителя весовой функции \(\rho\), многочлен \(Q\)
не может быть элементом ядра оператора \([T(0)]^{-1}T'\). Соответственно,
существует номер \(N\geqslant n\), для которого функциональная последовательность
\(\{Q_m\}_{m=0}^{\infty}\) вида
\[
	Q_m\rightleftharpoons\lambda_N^m\{[T(0)]^{-1}T'\}^mQ
\]
сойдётся в пространстве \(W_2^2[0,1]\) к нетривиальному кратному собственной
функции \(y_N\). При этом [\ref{prop:2.5}, \ref{prop:2.4}, \ref{prop:2.10}] число
нулей функции \(y_N\) не может превосходить числа знакоперемен многочлена \(Q\),
а тогда и величину \(n\). Объединяя сказанное, убеждаемся в выполнении равенства
\(N=n\) и наличии у собственной функции \(y_n\) в точности \(n\) нулей на интервале
\((0,1)\).

Распространение полученных результатов на общий случай \(\alpha\geqslant 0\)
проводится предельным переходом с учётом утверждений \ref{prop:2.4}
и \ref{prop:2.10}.
\end{proof}


\section{Спектральная периодичность и асимптотики собственных
значений}\label{par:3}
\subsection
Имеют место следующие два факта:

\subsubsection\label{prop:3.1}
{\itshape Пусть \(\{\lambda_n\}_{n=0}^{\infty}\) "--- последовательность
занумерованных в порядке возрастания собственных значений граничной задачи
\ref{par:1}.\ref{pt:1.1}\,\eqref{ekvac}, \ref{par:1}.\ref{pt:grusl}\,\eqref{grusl1}
при \(\alpha=0\), \(\beta=2/b\), а \(\{\mu_n\}_{n=0}^{\infty}\) "--- аналогичная
последовательность для граничной задачи того же типа при \(\alpha=0\),
\(\beta=2a/b\). Тогда независимо от выбора индекса \(n\in\mathbb N\) выполняется
равенство
\[
	\lambda_{\varkappa n}=(\varkappa/a^3)\,\mu_n.
\]
}

\begin{proof}
Ввиду очевидного выполнения искомого равенства для собственных значений
\(\lambda_0=\mu_0=0\), достаточно ограничиться рассмотрением случая \(n>0\).

Зафиксируем отвечающую собственному значению \(\mu_n>0\) собственную функцию
\(y\), имеющую на интервале \((0,1)\) в точности \(n\) различных нулей
и не обращающуюся в нуль на границе этого интервала [\ref{par:2}.\ref{prop:2.12},
\ref{par:2}.\ref{prop:2.4}]. Ввиду простоты собственного значения \(\mu_n\),
тождество \(P(x)\equiv 1-P(1-x)\) \cite[1.1]{VSh:2011} гарантирует,
что удовлетворяющая уравнению
\[
	[y'''-\mu_nPy]'+\mu_nPy'=0
\]
и граничным условиям \ref{par:1}.\ref{pt:grusl}\,\eqref{grusl1} собственная функция
\(y\) является относительно точки \(1/2\) либо чётной, либо нечётной.
Это наблюдение позволяет построить функцию \(z\in C^3[0,1]\), удовлетворяющую
следующим условиям:
\begin{enumerate}
\item При любом выборе индекса \(k\in\{0,\ldots,\varkappa-1\}\) функция \(z_k\)
вида
\[
	z_k(x)\rightleftharpoons z(k[a+b]+ax)
\]
совпадает с функцией \(y\) с точностью до знака.
\item При любом выборе индекса \(k\in\{1,\ldots,\varkappa-1\}\) на интервале
\((k[a+b]-b, k[a+b])\) выполняется тождество
\[
	|z(x)|\equiv \left|y(0)+\dfrac{y''(0)}{2a^2}\cdot(x-k[a+b]+b)\cdot
		(x-k[a+b])\right|.
\]
\end{enumerate}
Непосредственным вычислением с учётом факта самоподобия функции \(P\)
устанавливается, что функция \(z\) удовлетворяет уравнению
\[
	[z'''-(\varkappa/a^3)\,\mu_nPz]'+(\varkappa/a^3)\,\mu_nPz=0.
\]
Кроме того, из неравенства \(y(0)\cdot y''(0)<0\) [\ref{par:1}.\ref{pt:grusl}\,%
\eqref{grusl1}, \ref{par:2}.\ref{prop:1.1}] вытекает наличие у функции \(z\)
в точности \(\varkappa n\) нулей на интервале \((0,1)\). Тем самым, доказываемое
утверждение является верным [\ref{par:2}.\ref{prop:2.12}].
\end{proof}

\subsubsection\label{prop:3.2}
{\itshape Пусть \(\{\lambda_n\}_{n=0}^{\infty}\) "--- последовательность
занумерованных в порядке возрастания собственных значений граничной задачи
\ref{par:1}.\ref{pt:1.1}\,\eqref{ekvac}, \ref{par:1}.\ref{pt:grusl}\,\eqref{grusl1}
при \(\alpha=12/b^2\), \(\beta=6/b\), а \(\{\mu_n\}_{n=0}^{\infty}\) "---
аналогичная последовательность для граничной задачи того же типа
при \(\alpha=12a^2/b^2\), \(\beta=6a/b\). Тогда независимо от выбора индекса
\(n\in\mathbb N\) выполняется равенство
\[
	\lambda_{\varkappa (n+1)-1}=(\varkappa/a^3)\,\mu_n.
\]
}

\begin{proof}
Данное утверждение доказывается аналогичным утверждению \ref{prop:3.1} образом
с тем основным отличием, что при "`сшивке"' копий исходной собственной функции
используются не квадратичные, а кубические параболы видов
\begin{equation}\label{eq:parab}
	\zeta_k\cdot\left[\dfrac{y''(0)}{3a^2b}\cdot\left(\zeta_k^2-
		\dfrac{b^2}{4}\right)+\dfrac{2y(0)}{b}\right],
\end{equation}
где положено \(\zeta_k\rightleftharpoons x-k[a+b]+b/2\). Ввиду заведомого
различия знаков величин \(y(0)\) и \(y''(0)\) [\ref{par:1}.\ref{pt:grusl}\,%
\eqref{grusl1}, \ref{par:2}.\ref{prop:1.1}], каждая из парабол \eqref{eq:parab}
имеет на отвечающем ей интервале \((k[a+b]-b,k[a+b])\) единственный нуль.
Последнее означает наличие у функции \(z\) в точности \(\varkappa(n+1)-1\) нулей
на интервале \((0,1)\).
\end{proof}

\subsection
Имеет место следующий факт:

\subsubsection
{\itshape Пусть \(N:(0,+\infty)\to\mathbb N\) "--- считающая функция собственных
значений пучка \(T\). Тогда при \(\lambda\to+\infty\) справедливо асимптотическое
соотношение
\begin{equation}\label{eq:as1}
	N(\lambda)=\lambda^D\cdot[s(\ln\lambda)+o(1)],
\end{equation}
где \(D\rightleftharpoons\nu^{-1}\ln\varkappa\), \(\nu\rightleftharpoons
\ln\varkappa-3\ln a\), а \(s\) "--- \mbox{\(\nu\)-пе}\-ри\-оди\-че\-ская функция,
допускающая на периоде \([0,\nu]\) представление
\begin{equation}\label{eq:as2}
	s(t)\equiv e^{-Dt}\sigma(t),
\end{equation}
в котором \(\sigma\) "--- некоторая чисто сингулярная неубывающая функция.
}

\begin{proof}
Заметим [\ref{par:1}.\ref{pt:1.2}\,\eqref{kvf}], что замена значений параметров
\(\alpha\) и \(\beta\) приводит к возмущению операторов пучка \(T\) некоторым
оператором не превосходящего \(4\) ранга. Соответственно, главный член асимптотики
считающей функции \(N\) не зависит от выбора указанных значений. На протяжении
оставшейся части доказательства в качестве основной будет рассматриваться задача
вида \(\alpha=0\), \(\beta=2a/b\) с последовательностью собственных значений
\(\{\mu_n\}_{n=0}^{\infty}\). Последовательность собственных значений задачи
\(\alpha=0\), \(\beta=2/b\) при этом будет обозначаться через \(\{\lambda_n\}_{%
n=0}^{\infty}\).

Введём в рассмотрение последовательность заданных на отрезке \([0,\nu]\) функций
вида \(\sigma_k(t)\rightleftharpoons\varkappa^{-k}N(e^{k\nu+t})\). Заметим, что
независимо от выбора значений \(k,\,n\in\mathbb N\) и \(t\in [0,\nu]\) выполнение
неравенств
\[
	\mu_n<e^{k\nu+t}\leqslant\mu_{n+1}
\]
влечёт [\ref{prop:3.1}, \ref{par:1}.\ref{pt:1.2}\,\eqref{kvf}] выполнение
неравенств
\[
	\mu_{\varkappa n}\leqslant\lambda_{\varkappa n}<e^{(k+1)\nu+t}\leqslant
		\lambda_{\varkappa(n+1)}\leqslant\mu_{\varkappa(n+1)+2}.
\]
Таким образом, при любых \(k\in\mathbb N\) и \(t\in [0,\nu]\) выполняется
неравенство
\begin{equation}\label{eq:333}
	|\sigma_{k+1}(t)-\sigma_k(t)|\leqslant\varkappa^{-k},
\end{equation}
автоматически означающее равномерную сходимость последовательности
\(\{\sigma_k\}_{k=0}^{\infty}\) к некоторой функции \(\sigma\) со свойствами
\eqref{eq:as1}, \eqref{eq:as2}.

Далее, независимо от выбора значений \(k,\,n\in\mathbb N\) и \(t\in [0,\nu]\)
выполнение неравенств
\[
	\sup(\mu_{\varkappa (n+1)-1},\lambda_{\varkappa n})<e^{(k+1)\nu+t}
		\leqslant\mu_{\varkappa (n+1)}
\]
влечёт выполнение равенства \(\sigma_{k+1}(t)=\sigma_k(t)\). При этом, путём
почти дословного повторения рассуждений из доказательств утверждений
\cite[\S~5.1.1]{VSh:2011} и \cite[\S~5.2.1]{VSh:2011}, устанавливается
[\ref{prop:3.1}, \ref{prop:3.2}] ограниченность последовательностей частичных
сумм рядов
\begin{gather*}
	\sum\limits_{n=1}^{\infty}|\ln\mu_{\varkappa (n+1)-1}-
		\ln\mu_{\varkappa n}|,\\
	\sum\limits_{n=1}^{\infty}|\ln\lambda_{\varkappa n}-\ln\mu_{\varkappa n}|.
\end{gather*}
Соответственно, последовательность мер множеств вида
\[
	\{t\in[0,\nu]\::\:\sigma_{k+1}(t)\neq\sigma_k(t)\}
\]
имеет при \(k\to\infty\) асимптотику \(o(1)\), что влечёт [\eqref{eq:333}]
справедливость асимптотических соотношений
\begin{gather*}
	\|\sigma_{k+1}-\sigma_k\|_{L_2[0,1]}=o(\varkappa^{-k}),\\
	\|\sigma_k-\sigma\|_{L_2[0,1]}=o(\varkappa^{-k}).
\end{gather*}
Учёт признака сингулярности \cite[\S~4.1.3]{VSh:2011} и того обстоятельства,
что функции \(\sigma_k\) заведомо имеют не более \(O(\varkappa^k)\) точек разрыва
[\eqref{eq:as1}], завершает доказательство.
\end{proof}

\subsection
Таблицы из настоящего пункта содержат данные, относящиеся к уравнению, весовой
функцией в котором выступает обобщённая производная канторовой лестницы.
Данные таблицы \ref{tab:1} иллюстрируют утверждение \ref{prop:3.1}. Данные таблицы
\ref{tab:2} иллюстрируют утверждение \ref{prop:3.2}.

\begin{table}[t]
\begin{center}
\begin{tabular}{|r|r@{\;}c@{\;}l|r@{\;}c@{\;}l|r@{\;}c@{\;}l|}
\hline
{\(n\)}&\multicolumn{3}{|c|}{\(\mu_n\)}&
\multicolumn{3}{|c|}{\(54\mu_n\)}&\multicolumn{3}{|c|}{\(\lambda_n\)}\\ \hline
1&\(2,2131\cdot 10^1\)&\(\pm\)&\(10^{-3}\)&\(1,1951\cdot 10^3\)&\(\pm\)&
	\(10^{-1}\)&\(4,0965\cdot 10^1\)&\(\pm\)&\(10^{-3}\)\\
2&\(8,1717\cdot 10^2\)&\(\pm\)&\(10^{-2}\)&\(4,4127\cdot 10^4\)&\(\pm\)&\(10^0\)&
	\(1,1951\cdot 10^3\)&\(\pm\)&\(10^{-1}\)\\
3&\(3,175\cdot 10^3\)&\(\pm\)&\(10^0\)&\(1,714\cdot 10^5\)&\(\pm\)&\(10^2\)&
	\(3,867\cdot 10^3\)&\(\pm\)&\(10^0\)\\
4&\(3,849\cdot 10^4\)&\(\pm\)&\(10^1\)&\(2,078\cdot 10^6\)&\(\pm\)&\(10^3\)&
	\(4,412\cdot 10^4\)&\(\pm\)&\(10^1\)\\
\hline
\end{tabular}
\end{center}

\vspace{0.5cm}
\caption{Оценки первых собственных значений задач \(\alpha=0\), \(\beta=2\)
и \(\alpha=0\), \(\beta=6\) для случая \(\varkappa=2\), \(a=b=1/3\).}
\label{tab:1}
\end{table}

\begin{table}[t]
\begin{center}
\begin{tabular}{|r|r@{\;}c@{\;}l|r@{\;}c@{\;}l|r@{\;}c@{\;}l|}
\hline
{\(n\)}&\multicolumn{3}{|c|}{\(\mu_n\)}&
\multicolumn{3}{|c|}{\(54\mu_n\)}&\multicolumn{3}{|c|}{\(\lambda_n\)}\\ \hline
0&\(8,2987\cdot 10^0\)&\(\pm\)&\(10^{-4}\)&\(4,4813\cdot 10^2\)&\(\pm\)&
	\(10^{-2}\)&\(4,0965\cdot 10^1\)&\(\pm\)&\(10^{-3}\)\\
1&\(1,3784\cdot 10^2\)&\(\pm\)&\(10^{-2}\)&\(7,443\cdot 10^3\)&\(\pm\)&\(10^0\)&
	\(4,4813\cdot 10^2\)&\(\pm\)&\(10^{-2}\)\\
2&\(1,6311\cdot 10^3\)&\(\pm\)&\(10^{-1}\)&\(8,808\cdot 10^4\)&\(\pm\)&\(10^1\)&
	\(3,867\cdot 10^3\)&\(\pm\)&\(10^0\)\\
3&\(4,380\cdot 10^3\)&\(\pm\)&\(10^0\)&\(2,365\cdot 10^5\)&\(\pm\)&\(10^2\)&
	\(7,443\cdot 10^3\)&\(\pm\)&\(10^0\)\\
4&\(4,586\cdot 10^4\)&\(\pm\)&\(10^1\)&\(2,476\cdot 10^6\)&\(\pm\)&\(10^3\)&
	\(6,251\cdot 10^4\)&\(\pm\)&\(10^1\)\\
5&\(6,465\cdot 10^4\)&\(\pm\)&\(10^1\)&\(3,491\cdot 10^6\)&\(\pm\)&\(10^3\)&
	\(8,808\cdot 10^4\)&\(\pm\)&\(10^1\)\\
\hline
\end{tabular}
\end{center}

\vspace{0.5cm}
\caption{Оценки первых собственных значений задач \(\alpha=12\), \(\beta=6\)
и \(\alpha=108\), \(\beta=18\) для случая \(\varkappa=2\), \(a=b=1/3\).}
\label{tab:2}
\end{table}

\end{document}